\documentclass[12pt,leqno]{amsart}

\usepackage{amssymb,amsfonts,amsmath,amsthm}
\usepackage{latexsym}
\usepackage{verbatim}
\usepackage{epsfig}
\usepackage[active]{srcltx}
\usepackage[T1]{fontenc}
\usepackage[latin1]{inputenc}

\newtheorem{theorem}{Theorem}[section]

\newtheorem{proposition}[theorem]{Proposition}
\theoremstyle{definition}

\theoremstyle{remark}

\theoremstyle{remark}

\theoremstyle{remark}

\theoremstyle{remark}

\theoremstyle{remark}

\theoremstyle{remark}

\renewcommand{\Box}{\square}    
\renewcommand{\Bbb}{\mathbb}
\newcommand{\cal}{\mathcal}

\renewcommand{\int}{{\rm{int}}}

\newcommand{\Sing}{{\rm{Sing}}}

\newcommand{\id}{{\rm{id}}}

\newcommand{\Iso}{{\rm{Iso}}}
\newcommand{\im}{\mathop{\rm{im}}\nolimits}

\renewcommand{\ker}{\mathop{{\rm{ker}}}\nolimits}

\newcommand{\cl}{{\rm{closure}}}

\newcommand{\ity}{{\infty}}

\newcommand{\fin}{\hspace*{\fill}$\Box$}
\newcommand{\m}{\setminus}

\newcommand{\var}{{\rm{var}}}

\newcommand{\cA}{{\cal A}}

\newcommand{\cH}{{\cal H}}

\newcommand{\cS}{{\cal S}}

\newcommand{\cW}{{\cal W}}

\newcommand{\bC}{{\Bbb C}}

\newcommand{\bP}{{\Bbb P}}

\newcommand{\bZ}{{\Bbb Z}}
\newcommand{\bX}{{\Bbb X}}

\newcommand{\bH}{{\Bbb H}}

\begin{document}

\title[Minimal models of arrangements]{Complements of hypersurfaces, variation maps and minimal models of arrangements}
\author{Mihai Tib\u ar}

\address{Math\' ematiques, UMR 8524 CNRS,
Universit\'e des Sciences et Technologies de Lille, \  59655 Villeneuve d'Ascq, France.}
\email{tibar@agat.univ-lille1.fr}

\dedicatory{To Alexandru Dimca and \c Stefan Papadima, on the occasion of a great  anniversary}

 

\maketitle

\begin{abstract} 
We prove the minimality of the CW-complex structure for complements of hyperplane arrangements in $\bC^n$ by using the theory of Lefschetz pencils and results on the variation maps within a pencil of hyperplanes. This also provides a method to compute
the Betti numbers of complements of arrangements via global polar invariants. 
\end{abstract}

\setcounter{section}{0}
\section{Introduction}

 To study the topology of the complement $\bC^n \m V$ of an affine hypersurface $V\subset \bC^n$ one employs Morse theory, see for instance  Randell \cite{Ra}, or the Lefschetz method of scanning by pencils of hyperplanes, as done e.g. by Dimca and Papadima in \cite{DP}.  
Both methods yield in particular a CW-complex model of the complement $\bC^n \m V$. It was proved in the above two papers that whenever V is a union of hyperplanes, then there exists a CW-complex model which is {\em minimal}, in the sense that the number of $q$-cells equals the Betti number $b_q(\bC^n \m V)$, for any $q$. This notion of minimality was introduced by  Papadima and Suciu in \cite{PS} for studying the higher homotopy groups of complements of hyperplane arrangements.
We give here a new proof of the minimality by using another method. We first prove the following result:
 
 \begin{theorem} \label{t:main}
Let $\cA$ be an affine arrangement of hyperplanes, not necessarily central.
Let  $V_\cA \subset \bC^n$ denote the union of hyperplanes in $\cA$ and let $\cH$ be a generic hyperplane with respect to $\cA$. Then,   one has the isomorphisms of $\bZ$-modules: $H_j(\bC^n \m V_\cA) \simeq H_j(\cH \m V_\cA \cap \cH)$, for $j\le n-1$
and $H_n(\bC^n \m V_\cA) \simeq H_n(\bC^n \m V_\cA, \cH \m V_\cA \cap \cH)$.

Moreover, the complement $\bC^n \m V_\cA$ has a minimal model.
\end{theorem} 


  Our alternate proof uses the behaviour of the {\em variation maps} within a pencil of hyperplanes on $\bC^n\m V_\cA$. 
It is based on a particular case (see Theorem \ref{t:lef}) of a general result on vanishing cycles of pencils, which involves variation maps, proved in \cite{Ti-newton, Ti-top}. We discuss in \S \ref{hyp} some aspects of the topology of pencils on complements of affine hypersurfaces, extracted from  
a general theory of {\em non-generic Lefschetz pencils} of hypersurfaces, which we have developped in a series of papers \cite{Ti-newton, Ti-lef, Ti-top, Ti-surv,  Ti-book}.  
 In this context, we also give a method to compute inductively the betti numbers of complements of arrangements by using global polar invariants  \cite{Ti-imrn}. 

This question was brought to our attention by \c Stefan Papadima in spring 2000 in connexion with  \cite{PS}  (a preprint at that time) and with the earlier paper \cite{Ti-imrn} in which we construct CW-complex models for affine hypersurfaces  by using pencils and  global polar curves (see \S \ref{polar}). This note was essentially written in 2003  but not published ever since.  However, we think that it might be still of current interest also because of the recent proof by J. Huh \cite{Huh}   of a conjecture about the \emph{polar degree} stated by Dimca and Papadima \cite{DP}, of which one of the main ingredients is the non-generic Lefschetz pencil theory which we also use here.
  
\section{Complements of hypersurfaces}\label{hyp}    
Let $V = \{ f=0\}$ be a hypersurface in $\bC^n$. 
The complement $\bC^n \m V$ is a Stein manifold, since it can be viewed as the hypersurface $\{ t f(x) =1\}$ in $\bC^{n+1}$. It therefore has the homotopy type of a CW-complex of dimension $\le n$, by Hamm's result \cite{H1}.
For a generic hyperplane $\cH\in \bC^n$
 we have that the pair $(\bC^n \m V, \cH\m V\cap \cH)$ is $(n-1)$-connected, by Lefschetz type results \cite{H1,H2}, see also \cite[Thm.4.1]{Ti-surv}. One has the following well-known consequence:
  
\begin{proposition}\label{p:attach}
 The space $\bC^n \m V$ is obtained, up to homotopy type, by attaching to the slice $\cH\m V\cap \cH$ a certain number of $n$-cells.
 \fin
\end{proposition}

In general, hypersurface complements do not have minimal models:
  examples are given in \cite{PS}, one of the simplest being the case of
  the plane cusp $V= \{ x^2 - y^3 =0\}$. 
  
 It has been observed that the topology of the complement $\bC^n\m V$ depends on the singularities
of $V$ and also on their {\em position}, see \cite{Li1, Li2}. Moreover, 
if $V$ is not (stratified) transversal to the hyperplane at infinity,
then the non-transversality points may influence the topology, see \cite{Li1}, \cite{LT}.

\subsection{\ }
 A new viewpoint appeared more recently \cite{Ti-lef, LT}: consider a polynomial function $f: \bC^n\to \bC$ of which $V$ is a fiber, and relate
 the topology of the complement to the singularities of $f$.
  It is shown in \cite{LT} that one has two
situations: either $V$ is a general fiber of $f$ or a special one. For some fixed $f$, special (or ``atypical'') fibers are finitely many and have either singularities in $\bC^n$ or have, in some sense, singularities ``at infinity'' (see e.g. \cite{ST}). 
 We have:

\begin{proposition}\label{p:gen} \rm \cite{LT} \ \it
Let $V$ be a general fiber of some polynomial $f: \bC^n\to \bC$.
Then $\bC^n \m V$ is homotopy equivalent to the wedge $S^1 \vee S(V)$,
where $S(V)$ denotes the suspension over $V$. The cup-product in the cohomology ring of $\bC^n \m V$ is trivial.
\end{proposition}

 Even if in the above statement $V$ is non-singular, the complexity of the singularities at infinity of the polynomial $f$ influences the topology of $V$ (see \cite{ST}, \cite{LT} for examples).
   In certain situations, the general fiber of a polynomial function may be a bouquet of spheres of dimension $n-1$. It is the case when $f$ has isolated singularities at infinity. We send to \cite{Ti-surv} for a survey and more bibliography on singularities at infinity of polynomials. 
   
 When $V$ is an atypical fiber of a polynomial, we have the following result.

\begin{proposition}\label{p:aty} \rm \cite{LT} \ \it
Let $V = f^{-1}(0)$ be an atypical fiber of the polynomial function $f: \bC^{n+1} \to \bC$.
 If the general fibre of $f$ is 
$s$-connected, $s\ge 2$, then $\pi_i(\bC^{n+1}\setminus V) =0$, for $1<i\le 
s$, and $\pi_1(\bC^{n+1}\setminus V) =\bZ$.   
  \fin
\end{proposition} 
 
 In particular, if $f$ has ``isolated singularities at infinity'' then, the above discussion yields that $\pi_i(\bC^{n+1}\setminus V) =0$ for $1<i\le n-1$.

\section{Variation maps of pencils of affine hypersurfaces}\label{pencils}    

\subsection{Pencils with isolated singularities}\label{ss:pencils}
    
 The two methods of investigating the topology of complements, by Morse functions or by Lefschetz pencils, are actually close in spirit.
 The latter allows one to use the full power of complex geometry
 and we shall stick to it in this paper.
 
 In several recent papers we have introduced and used a general concept of {\em non-generic pencils of hypersurfaces} (e.g. \cite{Ti-compo, Ti-lef, Ti-surv}), which may have singularities in the axis. Here we only use pencils of hyperplanes and with
``no singularities in the axis'', as we describe in the following.   
 
 We consider our complement $\bC^n\m V$ as embedded into the projective space  $\bP^n$, and we identify it to $\bP^n\m (\bar V \cup H^\ity)$, where $\bar V$ denotes the projective closure of 
the affine hypersurface $V = f^{-1}(0)$ and $H^\ity$ is the hyperplane at infinity of $\bP^n$. Then consider the following pencil
of hyperplanes:
\begin{equation}\label{eq:pencil}
 l(x) - tx_0 =0, 
\end{equation}
where $t\in \bC$, $l : \bC^n \to \bC$ is a linear function and $x_0$ is the coordinate at infinity of $\bP^n$.   
  This pencil defines a holomorphic function $t:= l/x_0$ on $\bC^n = \bP^n\m H^\ity$, where $H^\ity$ denotes the hyperplane at infinity.
  Such a pencil is not generic with respect to the divisor $\bar V\cup H^\ity$ since the axis $A:= \{ l =x_0 = 0\}$ is included into 
 $H^\ity$ and hence $A$ is not
transversal to any Whitney stratification of the pair
$(\bP^n, \bar V \cup H^\ity)$. Nevertheless, we show that 
this pencil is {\em without singularities in the axis}, in the sense 
of \cite[Definitions 2.2, 2.3]{Ti-lef}.

 The projective hypersurface $\bar V \cup H^\ity$ has a canonical minimal Whitney stratification, which we denote by $\cW$. In particular, the intersection $\bar V \cap H^\ity$ is a union of strata. Then we consider the product stratification $\cW\times \bC$ in the product space $\bP^n \times \bC$. 
 
By a Bertini type result, there is a Zariski-open dense set $\Omega \subset \check \bP^{n-1}$ of linear forms $l : \bC^n \to \bC$ such that, for any $l\in \Omega$, the projective hyperplane $\{ l=0\} \subset H^\ity \simeq \bP^{n-1}$ is transversal within $H^\ity$
to all strata included into $\bar V \cap H^\ity$. In particular $\{ l=0\}$ avoids all point-strata inside $H^\ity$.

For such $l\in \Omega$, the hyperplane $\bH \subset \bP^n \times \bC$ defined by the equation (\ref{eq:pencil}) 
is transversal within  $\bP^n \times \bC$  to all product-strata
included into $H^\ity \times \bC$. Then the stratification $\cW\times \bC$ induces a stratification on $\bH$, call it $\cS$, which is also Whitney, by the transversality of the intersection.

  Moreover, $\cS$ has the property that all its strata which are included
  into $H^\ity \times \bC$ have a product structure, by the line $\bC$. It then follows that each member of the pencil
(i.e. for fixed $t \in \bC$) is transversal 
to all strata of $\cS$ included into $H^\ity \times \bC$. Equivalently, the projection to $\bC$ has no stratified singularities in the neighbourhood of $\bH \cap (H^\ity \times \bC)$.
 In such a case we say that the pencil  (\ref{eq:pencil})
 has no singularities in the axis.
It follows that this pencil can have singularities only outside the axis and that they are isolated. Namely, there are finitely many points on $V$ where the projection to the second factor $p: \bH \to \bC$ has a stratified singularity, with respect to the stratification $\cS$. 
The set of these points will be denoted by $\Sing_\cS p$.

\subsection{Variation maps} \label{var}


We recall from \cite{Ti-newton, Ti-top, Ti-book} and adapt to our case the construction of 
the global variation maps associated to a pencil.
Let us fix some notation. Let $X:= \bC^n \m V$ and note that $X$ can be identified to $\bH \cap ((\bC^n \m V)\times \bC)$. 

 For any $M\subset \bC$, we denote $\bH_M := p^{-1}(M)$ and $X_M := \bH_M \cap ((\bC^n \m V)\times \bC)$.
 Let $\Sing_\cS p = \cup_{i,j}\{a_{ij}\}$, where $\Lambda := p(\Sing_\cS p) = \{ a_1, \ldots , a_p\}$ and 
 $a_{ij}$ denotes some point of $\Sing_\cS p \cap p^{-1}(a_i)$.

  For $c\in \bC\m \Lambda$ we say that $\bH_c$, resp. $X_c$, is a {\em general fiber} of $p : \bH\to \bC$, resp. of $p_{|} : \bH \cap ((\bC^n \m V)\times \bC)\to \bC$. Indeed, $p_{|}$ can be identified to $l_| : \bC^n\m V \to \bC$ and $X_c$ is just $l^{-1}(c)\cap X$.

At some singularity $a_{ij}\in V$, we choose a ball $B_{ij}$ centered at $a_{ij}$. For a small enough radius of $B_{ij}$,
this is a ``Milnor ball" of the holomorphic function $p$ at $a_{ij}$.  Next we may take a small enough disc $D_i\subset \bC$ at $a_i \in \bC$, so that $(B_{ij}, D_i)$ is Milnor data for $p$ at $a_{ij}$. Moreover, we may do this for all (finitely many) singularities in the fiber $\bH_{a_i}$, keeping the same disc $D_i$, provided it is small enough.

Now the restriction of $p$ to $\bH_{D_i} \m \cup_j B_{ij}$ is a trivial fibration over $D_i$. One may construct a stratified vector field which trivializes this fibration and such that this vector field is tangent to the boundaries of the balls $\bH_{D_i} \cap \partial \bar B_{ij}$. Using this, we may also construct a geometric monodromy of the fibration $p_| : \bH_{\partial \bar D_i} \to \partial \bar D_i$ over the circle $\bar D_i$, such that this monodromy is the identity on the complement of the balls, $\bH_{\partial \bar D_i} \m \cup_j B_{ij}$. The same is then true, when replacing $\bH_{\partial \bar D_i}$ by $\bX_{\partial \bar D_i}$.

Fix some point $c_i \in \partial \bar D_i$. We have the geometric monodromy representation:
\[ \rho_i : \pi_1 (\partial \bar D_i, c_i) \to \Iso (X_{c_i}, X_{c_i}\m \cup_j B_{ij}), \]
where $\Iso (.,.)$ denotes the group of relative isotopy classes of stratified homeomorphisms (which are C$^\ity$ along each stratum).
It follows that the geometric monodromy restricted to $X_{c_i}\m \cup_j B_{ij}$ is the identity.

As shown above, we may identify the fiber $X_{c_i}\m \cup_j B_{ij}$ to the fiber $X_{a_i}\m \cup_j B_{ij}$ in the trivial fibration over $D_i$. Furthermore, in local coordinates at $a_{ij}$, $X_{a_i}$ is a germ of a complex analytic space; hence, for a small enough ball $B_{ij}$, the set $B_{ij}\cap X_{a_i}$ retracts to $\partial \bar B_{ij}\cap X_{a_i}$, by the local conical structure of analytic sets \cite{BV}.
  Therefore $X_{a_i}$ is homotopy equivalent, by retraction, to $X_{a_i}\m \cup_j B_{ij}$. 
  
\noindent
{\bf Notation}  \emph{Due to the above homotopy equivalences,  we shall freely use  $X_{a_i}^*$ as notation for $X_{c_i}\m \cup_j B_{ij}$ whenever we consider the pair $(X_{c_i}, X_{a_i})$.}

It then follows that the geometric monodromy induces an algebraic monodromy, in any dimension $q$:
\[ \nu_i \colon H_q (X_{c_i}, X_{a_i}^*; \bZ) \to  H_q (X_{c_i}, X_{a_i}^*; \bZ),\]
such that the restriction $\nu_i \colon H_q (X_{a_i}^*)\to  H_q(X_{a_i}^*)$ is the identity.
 
 Consequently, any relative cycle $\delta \in H_q (X_{c_i}, X_{a_i}^*; \bZ)$ is sent by the morphism $\nu_i - \id$ to an absolute cycle. In this way we define a {\em variation map}, for any $q\ge 0$:  
 
\begin{equation}\label{eq:var}
 \var_i : H_q(X_{c_i}, X_{a_i}^*; \bZ) \to H_q(X_{c_i}; \bZ).
\end{equation} 

Variation morphisms are basic ingredients in the description of the behaviour of vanishing cycles of global and local fibrations at singular fibers of holomorphic functions, see e.g. \cite{Mi}, \cite{La}, \cite{Si}, \cite[4.4]{Ti-compo}. Zariski already used $\nu_i - \id$ in dimension 2, in his well-known theorem for the fundamental group.
We shall use of \cite[Theorem 4.4]{Ti-top} in the following  form adapted to our particular case.
 
\begin{theorem}\label{t:lef}\rm \cite{Ti-newton, Ti-top} \it \ 
Let $V\subset \bC^n$, $l\in \Omega$ and let $X_c$ be a general member of the pencil, as above.
Then $H_q(X, X_{c})=0$ for $q\le n-1$ and the kernel of the surjection
$H_{n-1}(X_{c}) \twoheadrightarrow H_{n-1}(X)$
is generated by the images of the variation maps $\var_i$, for $i= \overline{1,p}$.
\end{theorem}

The first claim is also a consequence of the connectivity result stated in Proposition \ref{p:attach}. The second claim is highly nontrivial and is proved in \cite{Ti-top}. All the assumptions made in \cite[Theorem 4.4]{Ti-top} are clearly verified, except of one, which we still need to verify: $H_q(X_{c}, X_{a_i}^*)=0$ for $q\le n-2$.
This is indeed true by the following reason. In \cite[3.7, 3.9]{Ti-top}  it is shown that 
the named condition is satisfied whenever $H_q(X_{D}, X_{c})=0$ for $q\le n-1$,
where $D$ is a small enough disc centered at some value $a\in \Lambda$.
But the later condition is fulfilled by our \cite[Corollary 2.7]{Ti-compo}, which is based on Hamm and L\^e's results in \cite{HL}.

\subsection{Number of cells and polar invariants}\label{polar}

 Vanishing cycles in a pencil of hypersurfaces have been investigated
 in large generality, for example in \cite{ST, Ti-compo, Ti-imrn, Ti-surv}.  If the hypersurface $V\subset \bC^n$ is given by $f=0$ then, for some linear function $l$, one defines the global polar variety:
 \[ \Gamma (l,f):= \cl \{ \Sing(l,f) \setminus \Sing f\} \subset \bC^n. 
 \]
 
By the {\em global polar curve lemma} \cite[Lemma 2.4]{Ti-compo}, it follows that $\Gamma (l,f)$ is either empty or it is a curve, provided that $l$ is general enough.  This means that $l$ can be taken out of a Zariski-open set $\tilde \Omega\subset \Omega \subset\check \bP^{n-1}$, see {\em loc.cit.})  Global polar curves appeared for the first time in \cite{Ti-compo} in the study of the topology at infinity of polynomial functions. Local polar varieties have been introduced by L\^e D.T. and B. Teissier and are currently used in the literature. We refer the reader to \cite{Ti-imrn, ST-mon, Ti-surv} for different aspects of global polar curves.
 
By \cite[Theorem 4.6]{Ti-compo} and especially \cite[Corollary 4.3]{Ti-surv} we have that the Betti number $b_n(X, X_c)$ is equal to 
$\lambda := \sum_{i=1}^p \lambda_{a_i}$,
where $\lambda_{a_i}$ is the polar number at the atypical value of the pencil $a_i$.
According to \cite[Definition 3.5]{Ti-imrn}, $\lambda_{a_i}$ is a non-negative integer equal to the following difference of intersection multiplicities:
\[ \lambda_{a_i} = \int (\Gamma (l,f), X_{c_i}) - \int (\Gamma (l,f), X_{a_i}),\]
where $c_i$ is a nearby typical value of the pencil.

The difference of intersection numbers appears as follows. First observe that $\Gamma (l,f)$ does not intersect some small neighbourhood of 
$\bar V \cap H^\ity$. Next, the curve $\Gamma (l,f)$ is algebraic, therefore it intersects $V$ at a finite number of points.  It is a general fact proved by L\^e D.T \cite{Le} that these points are among  the stratified singularities of the restriction of the function $l$ to $V$. On the other hand, these singularities are isolated since $l$ is general. In the particular case of arrangements of hyperplanes, the stratified singularities of 
the restriction $l_{|V}$ are precisely the $0$-dimensional strata
of the canonical Whitney stratification of $V$. So the polar curve is
eventually non-empty in the neghbourhood of these points.
As $c_i$ tends to $a_i$, the points of intersection of  $\Gamma (l,f)$ with $X_{c_i}$, in some neghbourhood of some $0$-dimensional stratum of $V$ which is also on $\bar X_{a_i}$, tend to this point-stratum. Consequently, 
there is loss of intersection multiplicity from $\int (\Gamma (l,f), X_{c_i})$ to  $\int (\Gamma (l,f), X_{a_i})$ and this loss is localized near the point-strata of $V$.  

 Moreover, the space $X$ is obtained from the slice $X_c$ by attaching 
cells of dimension $n$ only, by \cite{Ti-surv}, \cite{ST-mero}, see also \cite[Theorem 9.3.1]{Ti-book}. We thus have a geometric interpretation  of the topological quotient space $X/X_c$ as a bouquet of $\lambda = b_n(X, X_c)$ $n$-spheres. By repeated slicing we get similar formulas in lower dimensions.
In case of complements of hyperplane arrangements, we shall in the next section that 
the relative betti number $b_n(X, X_c)$ equals the absolute betti number $b_n(X)$.

\section{Proof of Theorem \ref{t:main}}\label{s:arran}
We proceed by induction on the dimension.
Our arrangement of hyperplanes $\cA$ defines a natural Whitney stratification $\cW' = \{ W_B\}_{B\subset \cA}$ on $V_\cA$ which is also the coarsest one. More explicitly, the strata are defined as follows. Let $V_B$ denote the intersection of all hyperplanes corresponding to the indices of some subset $B \subset \cA$. Then $W_B := V_B \m \cup_{C\varsubsetneq B} V_C$. This stratification is Whitney since along any stratum $W_B$, by some analytic local change of coordinates, the space $V_\cA$ has the product structure $\{$transversal slice$\}\times W_B$. 

Since the hyperplane at infinity $H^\ity \subset \bP^n$ is transversal to all the strata, 
the induced natural stratification on $\bar V_\cA \cap H^\ity$
is Whitney and it is the coarsest one. This is what we have denoted by $\cW$ in \S \ref{ss:pencils}.

Let $l\in \Omega$ define a generic pencil of hyperplanes in $\bC^n$, as in \S \ref{ss:pencils}. We have seen before that the genericity of the pencil amounts to the condition that  the axis of the pencil $A = \{ l=0\} \cap H^\ity$ is transversal to all the strata of $\cW$.

Let $\cH$ denote a generic member of the pencil. 
  By Proposition \ref{p:attach}, we get that the long exact sequence of the pair $(\bC^n \m V_\cA, \cH \m V_\cA \cap \cH)$ splits into the isomorphisms $H_j(\bC^n \m V_\cA) \simeq H_j(\cH \m V_\cA \cap \cH)$, for $j\le n-1$, 
and the following exact sequence:
\begin{equation} \label{eq:long} 
 \begin{array}{l}0\to H_n(\bC^n \m V_\cA) \to H_n(\bC^n \m V_\cA, \cH \m V_\cA \cap \cH) \to \ \ \ \ \ \ \ \ \ \ \ \ \\  \ \ \ \ \ \ \ \ \ \ \ \
\to H_{n-1}(\cH \m V_\cA \cap \cH) \stackrel{\iota_*}{\to} H_{n-1}(\bC^n \m V_\cA)\to 0.\end{array} 
\end{equation} 
 We claim that $\iota_*$ is injective. By Theorem \ref{t:lef} we have that $\ker \iota_* = \sum_{i=1}^p \im(\var_i)$. In our case, we may show that  $\var_i$ is {\em trivial}, for any $i$. 
 Our pencil has no singularities in the axis, it is a pencil of hyperplanes and $V_\cA$ is a union of hyperplanes too. It follows that the singularities of the pencil are exacly the point-strata of the canonical stratification $\cW$ of $V_\cA$. Then the atypical members of the pencil are those which pass through such points. The pencil can be chosen generic enough such that each member of it contains at most one such point-stratum. 
 
  Let us focus on some atypical value $a_i$. We may assume, without affecting the generality, that the singularity of $\bar X_{a_i}$ is the origin of $\bC^n$. Consider the map germ $(\bC^n,0)\to (\bC^n,0)$ such that $x\mapsto x\exp(2i\pi t)$ for any coordinate $x$. Taking $t$ as parameter, this defines a family of diffeomorphisms which preserve the arrangement $V_\cA$ and its complement 
 $\bC^n \m V_\cA$, and moves the hyperplane $X_{c_i}$ of our pencil into the hyperplane  $X_{\exp(2i\pi t)c_i}$, over the circle $\partial \bar D_i \subset \bC$. 
   For $t=1$, this yields a geometric monodromy of $X_{c_i}$ around the value $a_i$, at the origin of $\bC^n$. 
 
 By its definition, this geometric monodromy is the identity on the hyperplane $X_{c_i}$ and therefore also on $X_{a_i}^* \subset X_{c_i}$ (see the definition of the notation $X_{a_i}^*$ at \S \ref{var}). 
  It then follows (from the definition of the variation map, see \ref{var})  that the variation of this monodromy is trivial, i.e. $\im (\var_i) = 0$. 
  We have proved in this way that $\ker \iota_* = 0$, which also means that the above exact sequence (\ref{eq:long}) splits in the middle. This proves the second part of our first statement. 
  
 By the attaching result discussed at the end of \S \ref{pencils}  (see also \cite[Theorem 9.3.1]{Ti-book})  we get that the number of the $n$-cells attached to $\cH \m V_\cA \cap \cH$ in order to obtain $\bC^n \m V_\cA$ is equal to $b_n(\bC^n \m V_\cA)$.
 The minimal model claim follows then by iterated slicing.


\end{document}